\documentclass[final,twoside,11pt]{entics}
\usepackage{enticsmacro}
\usepackage{graphicx}
\usepackage{amsmath}
\usepackage[all]{xy}
\usepackage{hyperref}
\newcommand{\ua}{\mathord{\uparrow}}
\newcommand{\da}{\mathord{\downarrow}}
\newcommand{\rom}[1]{\rm{\uppercase\expandafter{\romannumeral #1}}}

\def\conf{ISDT 2022} 	
\def\myconf{ISDT9}	
\volume{2}			
\def\volu{2}			
\def\papno{1}			
\def\mydoi{10354}
\def\copynames{X. Chu, Q. Li} 
\def\runtitle{The $d^{*}$-space}
\def\lastname{Chu and Li}

\begin{document}
\begin{frontmatter}
  \title{The $d^{*}$-space}
  \author{Xiangping Chu\thanksref{myemail}}	
   \author{Qingguo Li\thanksref{ALL}\thanksref{coemail}}		
   \address{School of Mathematics\\ Hunan University\\				
    Changsha, Hunan, China}  							
  \thanks[ALL]{This work is supported by the National Natural Science Foundation of China (No.~12231007)}   
   \thanks[myemail]{Email: \href{mailto:myuserid@mydept.myinst.myedu} {\texttt{\normalshape
        XP1997@hnu.edu.cn}}}
  \thanks[coemail]{Corresponding author, Email:  \href{mailto:couserid@codept.coinst.coedu} {\texttt{\normalshape
        liqingguoli@aliyun.com}}}
\begin{abstract}
 In this paper, we introduce the concept of $d^{\ast}$-spaces. We find that strong $d$-spaces are $d^{\ast}$-spaces, but the converse does not hold. We give a characterization for a topological space to be a $d^{\ast}$-space. We prove that each retract of a $d^{\ast}$-space is a $d^{\ast}$-space. We obtain the result that for any $T_{0}$ space $X$ and $Y$, if the function space $TOP(X,Y)$ endowed with the Isbell topology is a $d^{\ast}$-space, then $Y$ is a $d^{\ast}$-space. We also show that for any $T_{0}$ space $X$, if the Smyth power space $Q_{v}(X)$ is a $d^{\ast}$-space, then $X$ is a $d^{\ast}$-space. Finally, we give an example of a $d^{\ast}$-space $X$ whose Smyth power space $Q_{v}(X)$ is not a $d^{\ast}$-space, showing that the converse is false.
\end{abstract}
\end{frontmatter}
\section{Introduction}\label{intro}
Domain theory, initially developed by Dana Scott, has been exceedingly active at the interface between Mathematics and Theoretical Computer Science as the denotational semantics of functional programming language. Scott proved that a domain endowed with the Scott topology induced by the order of the domain is sober. As the sobriety can give rise to a categorical equivalence between topological spaces and certain frames, it is a very important topological property in domain theory and non-Hausdorff topological space. Other than sobriety, there were some other weaker properties put forward and investigated by some scholars, for example, $d$-space, well-filtered space, etc. In \cite{9}, Xu and Zhao introduced the concept of strong $d$-spaces lying between $d$-spaces and $T_{1}$ spaces. Strong $d$-spaces are d-spaces, but the converse is not true. Well-filtered spaces and strong $d$-spaces do not have inclusive relationships. However, a strong $d$-space is a well-filtered space when  equipped with the Scott topology. Moreover, every coherent well-filtered space is a strong $d$-space.

In \cite{8}, Lu and Li introduced the concept of weak well-filtered space which is strictly weaker than well-filtered space and showed that the Johnstone's example is weak well-filtered but not well-filtered. In this paper, we introduce the concept of $d^{\ast}$-space. We find that the $d^{\ast}$-space is strictly weaker than the strong $d$-space and show that the Johnstone's example is a $d^{\ast}$-space but not a strong $d$-space.  We give a characterization for a topological space to be a $d^{\ast}$-space. We prove that a retract of a $d^{\ast}$-space is a $d^{\ast}$-space. We obtain the result that for any $T_{0}$ spaces $X$ and $Y$, if the function space $TOP(X,Y)$ endowed with the Isbell topology is a $d^{\ast}$-space, then $Y$ is a $d^{\ast}$-space. We also show that for any $T_{0}$ space $X$, if the Smyth power space $Q_{v}(X)$  is a $d^{\ast}$-space, then $X$ is a $d^{\ast}$-space. Meanwhile, we give a counterexample to illustrate that conversely, for a $d^{\ast}$-space $X$, the Smyth power space $Q_{v}(X)$ may not be a $d^{\ast}$-space.


\section{Preliminaries}
We refer to [1,3] for the standard definitions and notation of order theory and domain theory, and to [2,4] for those of topologies.

Let $P$ be a poset and $A\subseteq P$. We denote $\ua A=\{x\in P \mid x\geq a \mbox{ for some }a\in A\}$ and $\da A=\{x\in P \mid x\leq a \mbox{ for some } a\in A\}$. For every $a\in P$, we denote $\ua\{a\}=\ua a=\{x\in P \mid x\geq a\}$ and $\da\{a\}=\da a=\{x\in P \mid x\leq a\}$. $A$ is called an \emph{upper set} (resp., a \emph{lower set}) if $A=\ua A$ (resp., $A=\da A$). $A$ is called \emph{directed} provided that it is nonempty and every finite subset of $A$ has an upper bound in $A$. The set of all directed sets of $P$ is denoted by $\mathcal{D}(P)$. Moreover, the set of all subsets of $P$ is denoted by $2^{P}$.

A poset $P$ is called a \emph{dcpo} if every directed subset $D$ in $P$ has a supremum. A subset $U$ of $P$ is called \emph{Scott open} if (1) $U=\ua U$ and (2) for any directed subset $D$ for which $\vee D$ exists, $\vee D\in U$ implies $D\cap U\neq \emptyset$. All Scott open subsets of $P$ form a topology, we call it the \emph{Scott topology} on $P$ and denoted by $\sigma(P)$. We denote $\Sigma P =(P,\sigma(P))$.

For a $T_{0}$ space $(X,\tau)$, let $\mathcal{O}(X)$ (resp., $\Gamma(X)$ ) be the set of all open subsets (resp., closed subsets) of $X$. For a subset $A$ of $X$, the closure of $A$ is denoted by $\mathrm{cl}_{\tau}(A)$ or $\overline{A}$.  We use $\leq_{\tau}$ to represent the \emph{specialization order} of $X$, that is, $x\leq_{\tau}y$ iff $x\in \overline{\{y\}}$. We denote $\da_{\tau}\{a\}=\da_{\tau}a=\{x\in P \mid x\leq_{\tau} a\}$. A nonempty subset $A$ of $X$ is \emph{irreducible} if $A\subseteq B\cup C$ for closed subsets $B$ and $C$ implies $A\subseteq B$ or $A\subseteq C$.  A subset $B$ of $X$ is called \emph{saturated} if $B$ equals the intersection of all open sets containing it (equivalently, $B$ is an upper set in the specialization order). For a topological space $X$, we denote the set of all nonempty compact saturated subsets of $X$ with the order reverse to containment, i.e., $K_{1}\leqslant K_{2}$ iff $K_{2}\subseteq K_{1}$  by $Q(X)$. We consider the \emph{upper Vietoris topology} $\upsilon$ on $Q(X)$  generated by the sets $\Box U=\{K \in Q(X) \mid K \subseteq U\}$, where $U$ ranges over the open subsets of $X$. We use $Q_{\upsilon}(X)$ to denote the resulting topological space. A $T_{0}$ space $X$ is called a \emph{$d$-space} (i.e., \emph{monotone convergence space}) if $X$ (with the specialization order) is a $dcpo$ and $\mathcal{O}(X)\subseteq \sigma(X)$.

For a $T_{0}$ space $X$, let $\mathcal{K}$ be a filtered family under the inclusion order in $Q(X)$, which is denoted by $\mathcal{K}\subseteq_{filt}Q(X)$, i.e., for any $K_{1},K_{2}\in Q(X)$, there exists $K_{3}\in Q(X)$ such that $K_{3}\subseteq K_{1}\cap K_{2}$. $X$ is called \emph{well-filtered} if for any open subset $U$ and any $\mathcal{K}\subseteq_{filt}Q(X)$, $\bigcap \mathcal{K}\subseteq U$ implies $K\subseteq U$ for some $K\in \mathcal{K}$. $X$ is called \emph {coherent} if the intersection of two compact saturated subset is compact.

\begin{definition}\cite{9}\label{}
A $T_{0}$ space $X$ is called a \emph{strong $d$-space} if for any $D\in \mathcal{D}(X)$, $x\in X$ and $U\in \mathcal{O}(X)$, $\bigcap\limits_{d\in D}{\uparrow}d\cap{\uparrow}x\subseteq U$ implies ${\uparrow}d \cap{\uparrow}x\subseteq U$ for some $d\in D$.
\end{definition}

\begin{definition}\cite{10}
 A poset $L$ is called a \emph{consistent dcpo} if for any directed subset $D$ of $L$ with $\bigcap\limits_{d\in D}{\uparrow}d\neq\varnothing$, $D$ has a least upper bound in $L$.
\end{definition}

\begin{definition}\cite{8}
 A topological space $(X,\tau)$ is called \emph{weak well-filtered} if, whenever a nonempty open set $U$ contains a filtered intersection $\bigcap\limits_{i\in I}Q_{i}$ of compact saturated subsets, then $U$ contains $Q_{i}$ for some $i \in I$.
\end{definition}

\begin{definition}\cite{4}
 A retract of a topological space $Y$ is a topological space $X$ such that there are two continuous maps $s: X \rightarrow Y$ and $r : Y \rightarrow X$ such that $r \circ s=\operatorname{id}_{X}$.
\end{definition}
\begin{theorem} \cite{5}
Let $X$ be a topological space and $\mathcal{A}$ an irreducible subset of the Smyth power space $Q_{\upsilon}(X)$. Then every closed set $C\subseteq X$ that meets all members of $\mathcal{A}$ contains a minimal irreducible closed subset $A$ that meets all members of $\mathcal{A}$.
\end{theorem}
\begin{theorem} \cite{12}
 A topological space $X$ is well-filtered iff its upper space $Q(X)$ is well-filtered.
\end{theorem}

\section{Main result}
In this section, we introduce the $d^{\ast}$-space inspired by Lu and Li's work\cite{8} on weak well-filtered space.
\begin{definition}
A $T_{0}$ space $X$ is called a \emph{$d^{\ast}$-space} if for any $D\in\mathcal{D}(X)$, $x\in X$, $U\in \mathcal{O}(X)\backslash\{\varnothing\}$, $\bigcap\limits_{d\in D}{\uparrow}d\cap{\uparrow}x\subseteq U$ implies ${\uparrow}d \cap{\uparrow}x\subseteq U$ for some $d\in D$.
\end{definition}

\begin{remark}\label{1}
A $T_{0}$ space $X$ is a $d^{\ast}$-space iff for any $D\in \mathcal{D}(X)$, $x\in X$, $A\in \Gamma(X)\backslash\{X\}$, if ${\uparrow}d \cap{\uparrow}x\cap A\neq\varnothing$ for all $d \in D$, then $\bigcap\limits_{d\in D}{\uparrow}d\cap{\uparrow}x\cap A\neq\varnothing$.
Obviously, every coherent weak well-filtered space is a $d^{\ast}$-space.
\end{remark}

Clearly every strong $d$-space is a $d^{\ast}$-space. The following example shows that a $d^{\ast}$-space may not be a strong $d$-space.
\begin{example}\label{7}
Consider the natural number $N$ with the usual order . It is easy to verify that $(N, \sigma(N))$ is a $d^{\ast}$-space. Obviously, $(N, \sigma(N))$ is not a $d$-space. Hence, it is not a strong $d$-space.
\end{example}

\begin{proposition}
Let $(X,\tau)$ be a $d^{\ast}$-space. Then $\Omega(X)$ is a consistent dcpo and $\tau\subseteq\sigma(\Omega(X))$, where $\Omega(X)=(X,\leq_{\tau})$.
\end{proposition}
\begin{proof}
Let $D$ be a directed subset of $X$ and $\bigcap\limits_{d\in D}{\uparrow_{\tau}}d\neq\varnothing$. Suppose that $\sup D$ does not exist. Then for all ${x{\in}\bigcap\limits_{d\in D}{\uparrow}_{\tau}d}$, there exists ${y{\in}\bigcap\limits_{d\in D}{\uparrow_{\tau}}d}$ such that $x\nleq y$, it follows that $x\in X\backslash{\downarrow}_{\tau}y$ and $D\subseteq{\downarrow_{\tau}}y$. Thus $\mbox{cl}_{\tau}{(D)}\subseteq{\downarrow_{\tau}}y$ and $x\notin \mbox{cl}_{\tau}{(D)}$. So $\bigcap\limits_{d\in D}{\uparrow}_{\tau} d\cap\mbox {cl}_{\tau}(D)=\varnothing$. Since $D$ is directed, it is nonempty, there exists $d_{0}\in D$. Note that $\bigcap\limits_{d\in D}{\uparrow}_{\tau}d\cap{\uparrow}_{\tau}d_{0}\cap \mbox {cl}_{\tau}(D)=\varnothing$. Obviously, $\mbox {cl}_{\tau}{(D)}\neq X$. Therefore, there exists $d\in D$ such that ${\uparrow_{\tau}}d\cap{\uparrow_{\tau}} d_{0}\cap\mbox{cl}_{\tau}{(D)}=\varnothing$ by Remark \ref{1}. Again by the fact that $D$ is directed, there exists $d^{\prime}\in D$ such that $d\leq d^{\prime}$ and $d_{0}\leq d^{\prime}$. Thus $d^{\prime}\notin\mbox{cl}_{\tau}{(D)}$, which is a contradiction. Hence, $\Omega(X)$ is a consistent dcpo.

Suppose that $U\in\tau\backslash\ \{\varnothing\}$. Obviously, it is an upper set in the specialization order. Let $D$ be a directed set on $\Omega(X)$ with $\sup D \in U$. Then $\bigcap\limits_{d\in D}\uparrow_{\tau}d= \uparrow_{\tau} \sup D \subseteq U$. Thus there exists $d\in D$ such that $\uparrow_{\tau} d\subseteq U$, that is $D\cap U\neq \varnothing$. Therefore, $U\in\sigma(\Omega(X))$.
\end{proof}
\begin{example}
 Let us consider $P=N\cup\{a\}$, where $N$ is the set of natural numbers with the usual order and for all $n\in N$, $n$ and $a$ are incomparable. It is clear that $P$ is a consistent dcpo. But $(P,\sigma(P))$ is not a $d^{\ast}$-space. Indeed, $\bigcap\limits_{n\in N} {\uparrow n}\subseteq\{a\}$, but for all $n\in N$, ${\uparrow n}\nsubseteq \{a\}$. It illustrated that a consistent dcpo with Scott topology may be not a $d^{\ast}$-space.
 \end{example}
\begin{theorem}\label{10}
For a {\rm dcpo} $P$, the following two conditions are equivalent.
\begin{enumerate}
\item $\Sigma P$ is a $d^{\ast}$-space.
\item For any $A\in {\Gamma}(P)\backslash\{P\}$ and $x\in P$, ${\downarrow}({\uparrow}x \cap A)\in \Gamma(P)$.
\end{enumerate}
\end{theorem}

\begin{proof}
$(1)\Rightarrow(2)$: Obviously, ${\downarrow}({\uparrow} x\cap A)$ is a lower set. We only need to prove that $\sup D\in{\downarrow}({\uparrow}x \cap A)$ for any directed subset $D\subseteq{\downarrow}({\uparrow}x \cap A)$ and $\sup D$ exists.  Suppose not, we have that ${\uparrow}\sup D\cap{\uparrow}x \cap A=\varnothing$, that is $\bigcap\limits_{d\in D}{\uparrow} d\cap{\uparrow}x \cap A=\varnothing$. Since $A\neq P$, there exists $d_{0}\in D$ such that ${\uparrow}d_{0}\cap{\uparrow}x \cap A=\varnothing$ by Remark \ref{1}. So $d_{0}\notin{\downarrow}({\uparrow}x\cap A)$, which is a contradiction.

$(2)\Rightarrow(1)$: Let $D$ be a directed subset of $P$ and $\bigcap\limits_{d\in D}{\uparrow}d\cap{\uparrow}x\subseteq U$ for any $x\in P$ and any nonempty open subset $U$ of $P$. Assume that for every $d\in D$, ${\uparrow}d \cap{\uparrow}x\nsubseteq U$. Then ${\uparrow}d \cap{\uparrow}x\cap P\backslash\ U\neq\varnothing$. Thus $d\in{\downarrow}({\uparrow}x\cap P\backslash\ U)$, this implies that $D\subseteq {\downarrow}({\uparrow}x\cap P\backslash\ U)$. Since $P$ is a $dcpo$,  $\sup D$ exists. By~(2), we have $\sup D\in{\downarrow}({\uparrow}x\cap P\backslash\ U)$, that is ${\uparrow}\sup D\cap{\uparrow}x\cap P\backslash\ U\neq\varnothing$. Therefore, $\bigcap\limits_{d\in D}{\uparrow} d\cap{\uparrow}x \cap P\backslash\ U\neq\varnothing$, which is a contradiction. So $\Sigma P$ is a $d^{\ast}$-space.
\end{proof}

\begin{proposition}\label{11}
Let $P$ be a \rm{dcpo}. For any $A\in{\Gamma}(P)\backslash\{P\}$ and any $K\in Q(P)$, ${\downarrow}(K\cap A)\in\Gamma( P)$. Then $(P,\sigma(P))$ is weak well-filtered.
\end{proposition}
\begin{proof}
Let $\mathcal{K}$ be a filtered compact saturated subset family of $P$ and $\bigcap\mathcal{K}\subseteq U$ for any nonempty open subset $U$. We need to prove that there exists a compact saturated subset $K$ in $\mathcal{K}$ such that $K\subseteq U$. Suppose not, for any $K\in\mathcal{K}$, $K\nsubseteq U$. Then $K\cap P\backslash U\neq\varnothing$. By Rudin Lemma, there exists a minimal closed subset $C\subseteq P\backslash U$ such that for all $K\in\mathcal{K}$, $K\cap C \neq\varnothing$. For any $K^{\prime}\in\mathcal{K}$, there exists $K^{\prime\prime}\in\mathcal{K}$ such that $K^{\prime\prime}\subseteq K\cap K^{\prime}$. Then $\varnothing\neq C\cap K^{\prime\prime}\subseteq C\cap K\cap K^{\prime}\subseteq{\downarrow} (C\cap K)\cap K^{\prime}$. Thus ${\downarrow} (C\cap K)\cap K^{\prime}\neq\varnothing$. Since $C$ is a minimal closed subset, we have ${\downarrow} (C\cap K)=C$. So there exists a maximal element $x_{0}$ in $C$ because $P$ is a $dcpo$ and $C\in \Gamma(P)$. Note that $x_{0}\in C={\downarrow}(C\cap K)$. Hence, for all $K\in\mathcal{K}$,  there exists $a_{K}\in C\cap K $ such that $x_{0}\leq a_{K}$. By the fact that $x_{0}$ is a  maximal element, we know $x_{0}=a_{K}$. So $x_{0}\in\bigcap\mathcal{K}\subseteq U\subseteq X\backslash C$, which is a contradiction.
\end{proof}
\begin{lemma}\cite{11}\label{12}
For a poset $P$ and $A\in\Gamma(P)$, the following two conditions are equivalent.
\begin{enumerate}
\item${\downarrow}({\uparrow} x\cap A)\in\Gamma( P)$ for all $x\in P$.
\item${\downarrow}(K\cap A)=\bigcup\limits_{k\in K}{\downarrow}({\uparrow}k\cap A)\in\Gamma(P)$ for all $K\in Q(\Sigma P)$.
\end{enumerate}
\end{lemma}

\begin{corollary}\label{13}
Let $L$ be a \rm{dcpo}. If $(L,\sigma(L))$ is a $d^{*}$-space, then $(L,\sigma(L))$ is weak well-filtered.
\end{corollary}

\begin{proof}
Suppose $A\in{\Gamma}(P)\backslash\{P\}$ and $K\in Q(P)$. We need to show that ${\downarrow}(K\cap A)$ in $\Gamma( P)$. Since $(L,\sigma(L))$ is a $d^{*}$-space, by Theorem \ref{10}, for any $A\in {\Gamma}(P)\backslash\{P\}$ and $x\in P$(especially $x$ in $K$), ${\downarrow}({\uparrow}x \cap A)\in \Gamma(P)$. By Lemma \ref{11}, ${\downarrow}(K\cap A)\in \Gamma(P)$. By Proposition \ref{12}, $(L,\sigma(L))$ is weak well-filtered.
\end{proof}
\begin{problem}
In Corollary~\ref{13}, if dcpo is replaced by poset, does the conclusion still hold?
\end{problem}

From the example below, we can see that Johnstone dcpo endowed with the Scott topology is a $d^{*}$-space, but it is not a strong $d$-space.
\begin{example}
Recall the dcpo constructed by Johnstone in \cite{4}, which is defined as $\mathbb{J}=N\times(N\cup\{\infty\})$, with the order defined by $(j,k)\leq (m,n)$ iff $j=m$ and $k\leq n$ or $n=\infty$ and $k\leq m$. Obviously, the Johnstone space $\Sigma \mathbb{J}$ is a $d$-space. Clearly, $\bigcap\limits_{n\geq 2}{\uparrow}(1,n)\cap{\uparrow}(2,2)=\varnothing$
but ${\uparrow}(1,n)\cap{\uparrow}(2, 2)=\{(m, \infty)\mid n \leq m\}\neq\varnothing$ for any $n\geq 2$. Hence, $(J, \sigma(J))$ is not a strong~$d$-space. Consider any directed subset $D=\{(x_{i},y_{i})\}_{i\in I}$, $t=(x,y)\in J$ and $U\in \sigma(J)\backslash\{\varnothing\}$. If $\bigcap\limits_{i\in I}{\uparrow}(x_{i},y_{i})\cap{\uparrow}(x,y)\subseteq U$, then we have $\bigcap\limits_{i\in I}{\uparrow}(x_{i},y_{i})\cap{\uparrow}(x,y)={\uparrow}\sup D\cap{\uparrow}t\subseteq U$. We need to consider the following two cases.\\
Case 1. If $\sup D \in D$, then there exists $i_{0}\in I$ such that $(x_{i_{0}},y_{i_{0}})=\sup D$; hence, ${\uparrow}(x_{i_{0}},y_{i_{0}})\cap{\uparrow} (x,y)\subseteq U$.\\
Case 2. If $\sup D \notin D$, then $D$ is a chain. Hence, there exists $n_{0}\in N$ such that $x_{i}=n_{0}$ for all $i\in I$. So $\sup D=(n_{0},\infty)$. For $t=(x,y)$, if $x=n_{0}$, then ${\uparrow}\sup D\cap{\uparrow}t={\uparrow}\sup D\subseteq U$, that is, $\sup D\in U$. Since $U$ is a Scott open subset, there exists a $i_{1}\in I$ such that $(n_{0},y_{i_{1}})\in U$. Hence, ${\uparrow}(n_{0},y_{i_{1}})\cap{\uparrow}(x,y)\subseteq {\uparrow}(n_{0},y_{i_{1}})\subseteq U$. If $x\neq n_{0}$, then we can find a $k_{0}\in N$ such that $(k _{0},\infty)\in U$ because $U$ is nonempty. Clearly, there exists $m_{0}\in N$ such that $(k_{0},m_{0})\in U$. So there exists $(l_{0},\infty)\in U$ where $m_{0}\leq l_{0}$ and $y\leq l_{0}$. Thus there exists $i\in I$ such that $l_{0}\leq y_{i}$. Therefore, ${\uparrow}(x_{i},y_{i})\cap{\uparrow}(x,y)=\{(b, \infty)\mid b\geq y_{i}$  and  $b\geq y \}\subseteq{\uparrow}(k_{0},m_{0})\subseteq U$.\\ So the Johnstone space $\Sigma J$ is a $d^{*}$-spaces.\\
\end{example}

\begin{theorem}
If $A$ is a saturated subspace of $d^{*}$-space $X$, then $A$ is a $d^{*}$-space.
\end{theorem}
\begin{proof}
 Let $D$ be a directed subset of $A$ and $\bigcap\limits_{d\in D}{\uparrow}_{A}d\cap{\uparrow}_{A}x\subseteq U$ for any $x\in A$ and any nonempty open set $U$ on $A$. Since $A$ is a saturated subspace of $d^{*}$-space $X$, there exists $V\in \mathcal{O}(X)$ such that $U=V\cap A$. We claim that ${\uparrow}_{A}x={\uparrow}x$. Clearly, ${\uparrow}_{A}x \subseteq {\uparrow}x$. Assume $y\in{\uparrow}x$. Then $x\leq y$. Since $A$ is saturated, we have $y\in A$. Hence, $y\in{\uparrow}_{A}x$. Thus $\bigcap\limits_{d\in D}{\uparrow}d\cap{\uparrow}x\subseteq V\cap A\subseteq V$ and $V\neq\varnothing$. So there exists $d\in D$ such that ${\uparrow}d\cap{\uparrow}x\subseteq V$. We conclude that ${\uparrow}_{A} d\cap{\uparrow}_{A}x\subseteq V\cap A=U$.
 \end{proof}
 \begin{theorem}
If $A$ is a closed subspace of $d^{*}$-space $X$, then $A$ is a $d^{*}$-space.
\end{theorem}
\begin{proof}
Let $D$ be a directed subset of $A$ and $\bigcap\limits_{d\in D}{\uparrow}_{A}d\cap{\uparrow}_{A}x\subseteq U$ for any $x\in A$ and any nonempty open set on $A$. Since $A$ is a closed subspace of $d^{*}$-space $X$, there exists $V\in \mathcal{O}(X)$ such that $U=V\cap A$. Thus $\bigcap\limits_{d\in D}{\uparrow}_{A} d \cap{\uparrow}_{A} x\subseteq V\cap A$. So $\bigcap\limits_{d\in D}{\uparrow}d\cap{\uparrow} x\subseteq V\cup (X\backslash A)$. Obviously, $V\cup(X\backslash A)\neq\varnothing$. Hence, there exists $d\in D$ such that ${\uparrow} d\cap{\uparrow} x\subseteq V\cup(X\backslash A)$. It follows that ${\uparrow}_{A}d\cap{\uparrow}_{A}x\subseteq V\cap A=U$.
\end{proof}
In \cite{6}, it has been verified that the product of two strong $d$-spaces need not be a strong $d$-space. The following example reveals that the product of two $d^{*}$-spaces need not be a $d^{*}$-space, either.
\begin{example}\label{6}
 We know that $\Sigma N $ and the Sierpinski space $\Sigma2$ are $d^{*}$-spaces, where $2$ is defined as $\{0, 1\}$ with the ordering $0<1$. We consider the product space $\Sigma N \times \Sigma2$. Let $\mathcal{D}=\{(n, 0)\mid n\in N\}$ and $x=(0,0)$, $U={\uparrow}(0,1)=N\times\{1\}$. Obviously, $\varnothing=\bigcap\limits_{n\in N}{\uparrow}(n,0)\cap{\uparrow}(0,0)\subseteq N\times\{1\}$. But for any $n\in N$, ${\uparrow}(n,0)\cap{\uparrow}(0,0)={\uparrow}(n,0)\nsubseteq N\times\{1\}$. Thus $\Sigma N \times\Sigma2$ is not a $d^{*}$-space.
\end{example}

\begin{theorem}
A retract of a $d^{*}$-space is  a $d^{*}$-space.
\end{theorem}
\begin{proof}
Suppose $X$ is a $d^{*}$-space and $Y$ is a retract of $X$. Then there are two continuous maps $f: X \rightarrow Y $ and $g : Y \rightarrow X $ such that $f \circ g=\operatorname {id}_{Y}$. Suppose $D$ is a directed subset of $Y$ with $\bigcap\limits_{d\in D}{\uparrow} d\cap{\uparrow}x\subseteq U$ for any $x\in Y$ and  any nonempty open set $U$ on $Y$. We can prove that $\bigcap\limits_{d\in D}{\uparrow}g(d)\cap{\uparrow}g(x)\subseteq f^{-1}(U)$. Indeed, assume $t\in\bigcap\limits_{d\in D}{\uparrow} g(d)\cap{\uparrow} g(x)$. Then for any $d\in D$, $g(d)\leq t$ and $g(x)\leq t$. So $d=f(g(d))\leq f(t)$ and $x=f(g(x))\leq f(t)$. Thus $f(t)\in\bigcap\limits_{d\in D}{\uparrow}d\cap{\uparrow}x$, it follows that $f(t)\in U$. So $t\in f^{-1}(U))$. Since $f$ is continuous and surjective, we have $f^{-1}(U)\in \mathcal{O}(X)$ and $f^{-1}(U)\neq\varnothing$. So there exists $ d\in D$ such that ${\uparrow}g(d)\cap{\uparrow}g(x)\subseteq f^{-1}(U)$ by the assumption that $X$ is a $d^{*}$-space. Assume $m\in{\uparrow}d \cap{\uparrow}x$. Since $g$ is monotone, we have $g(d)\leq g(m)$ and $g(x)\leq g(m)$. It follows that $g(m)\in f^{-1}(U)$. Thus $m=fg(m)\in U$. We conclude that ${\uparrow}d \cap{\uparrow}x\subseteq U$, that is  $Y$ is a $d^{*}$-space.
\end{proof}
Given topological spaces $X$ and $Y$, let ${TOP}(X, Y)$ be the set of all continuous functions from
$X$ to $Y$. The Isbell topology on the set ${TOP}(X, Y)$ is generated by the subsets of the form
$N(H\leftarrow V)=\{f\in {TOP}(X, Y)$ $:f ^{-1}(V)\in H\}$, where $H$ is a Scott open subset of the complete lattice $\mathcal{O}(X)$ and $V$ is open in the topological space $Y$. Let $[X, Y]$ denote $ {TOP}(X, Y)$ endowed with the Isbell topology.\\

\begin{lemma}\cite{7}\label{2}
For each $y \in Y$,consider the mapping $\xi : Y \rightarrow [X, Y]$ by $\xi (y)=\xi_y$, where $\xi_{y}(x) = y $ for all $x \in X$. Then $\xi$ is continuous.
\end{lemma}

\begin{lemma}\label{3}
Consider the mapping $F:[X, Y]\rightarrow Y$ by $F(f)=f(x_{0})$ for some fixed $x_{0}\in X$. Then $F(f)$ is continuous.
\end{lemma}
\begin{proof}
Suppose $U\in \mathcal{O}(Y)$. Then $F^{-1}(U)=\{f\in [X, Y]|F(f)\in U\}=\{f\in [X, Y]|f(x_{0})\in U\}=\{f\in [X, Y]|f^{-1}( U)\in\mathcal{N}(x_{0})\}$. Since $f$ is continuous, we have $f^{-1}( U)\in \mathcal{O}(X)$. Obviously, $\mathcal{N}(x_{0})$ the set of all open neighbourhood of $x_{0}$ is a Scott open subset of the complete lattice $\mathcal{O}(X)$. So $\{f\in [X, Y]|f^{-1}( U)\in\mathcal{N}(x_{0})\}$ is open in $[X, Y]$. Therefore, $F(f)$ is continuous.
\end{proof}

\begin{proposition}
For any ${T}_{0}$ topological spaces $X$ and $Y$, $Y$ is a retract of $[X, Y]$.
\end{proposition}
\begin{proof}
By Lemma \ref{2}  and Lemma \ref{3}, we can define $F:[X, Y]\rightarrow Y$ by $F(f)=f(x_{0})$ for some fixed $x_{0}\in X$, $\xi : Y \rightarrow [X, Y]$ by $\xi (y)=\xi _{y}$, where $\xi_{y}(x) = y $ for all $x \in X$. Obviously, they are continuous and $F\circ\xi(y)=F(\xi(y))=F(\xi_{y})=\xi_{y}(x_{0})=y$. We conclude that $F\circ\xi=\operatorname{id}_{Y}$. So $Y$ is a retract of $[X, Y]$.
\end{proof}
Now we make the conclusion below.
\begin{corollary}\label{4}
If  $[X, Y]$ is  a $d^{*}$-space, then $Y$ is a $d^{*}$-space.
\end{corollary}
We give the following example to illustrate that the converse of Corollary \ref{4} does not hold.
\begin{example}
The Sierpinski space $\Sigma2$ (Example \ref{6}) and $(N,\sigma(N))$ are $d^{*}$-spaces (Example \ref{7}). For  any $n\in N$, We define the function $f_{n} :\Sigma 2\rightarrow\Sigma N$ by

\[f_{n}(x)=
\begin{cases}
n, x=1;\\
0, \mbox{else}.
\end{cases}
\]
Obviously, $\{f_{n}\}_{n\in N}$ is directed in $TOP(\Sigma 2,\Sigma N)$. It is also clear that $\{\{0,1\}\}$ is a Scott open subset of the complete lattice $\mathcal{O}(2)$ and ${\uparrow}3$ is open in $\Sigma N$. Let $f=f_{0}$ and $U=N(\{\{0,1\}\} \leftarrow{\uparrow} 3)=\{f\in [\Sigma2\rightarrow\Sigma N]\mid f^{-1}({\uparrow}3)\in\{\{0,1\}\}\}=\{f\in [\Sigma2\rightarrow\Sigma N]\mid f^{-1}({\uparrow} 3)=\{0,1\}\}$. Then  $\varnothing=\bigcap\limits_{n\in N}{\uparrow}f_{n}\cap{\uparrow}f=\bigcap\limits_{n\in N}{\uparrow}f_{n}\subseteq U$. But for any $n\in N$,
\[
f_{n}^{-1}({\uparrow} 3)=
\begin{cases}
\varnothing,3>n;\\
\{1\},3\leq n.\\
\end{cases}
\]
 Thus for any $n\in N$, $f_{n}\notin U$. Hence, ${\uparrow} f_{n}\cap{\uparrow}f={\uparrow}f_{n}\nsubseteq U$. So $[\Sigma2\rightarrow\Sigma N]$ is not a $d^{*}$-space.
\end{example}
We know that if a topological space $X$ is well-filtered, then $Q(X)$ is a dcpo. However, the converse of this result does not hold.
\begin{example}
Let $X=N\cup\{\omega_{1},\omega_{2}\}$, where the natural number $N$ with the usual order and for all $n\in N$, $n<\omega_{1},\omega_{2}$, $\omega_{1}$ and $\omega_{2}$ are incomparable. Then $(X,\sigma)$ is a $T_{0}$ topological space. Obviously, $Q(X)=\{{\uparrow}x \mid x\in X\}\cup\{\{\omega_{1},\omega_{2}\}\}$, and $Q(X)$ is a dcpo. But $(X,\sigma)$ is not well-filtered because $X$ is not a dcpo.
\begin{figure*}[h]
  \centering
  \includegraphics[width=15cm]{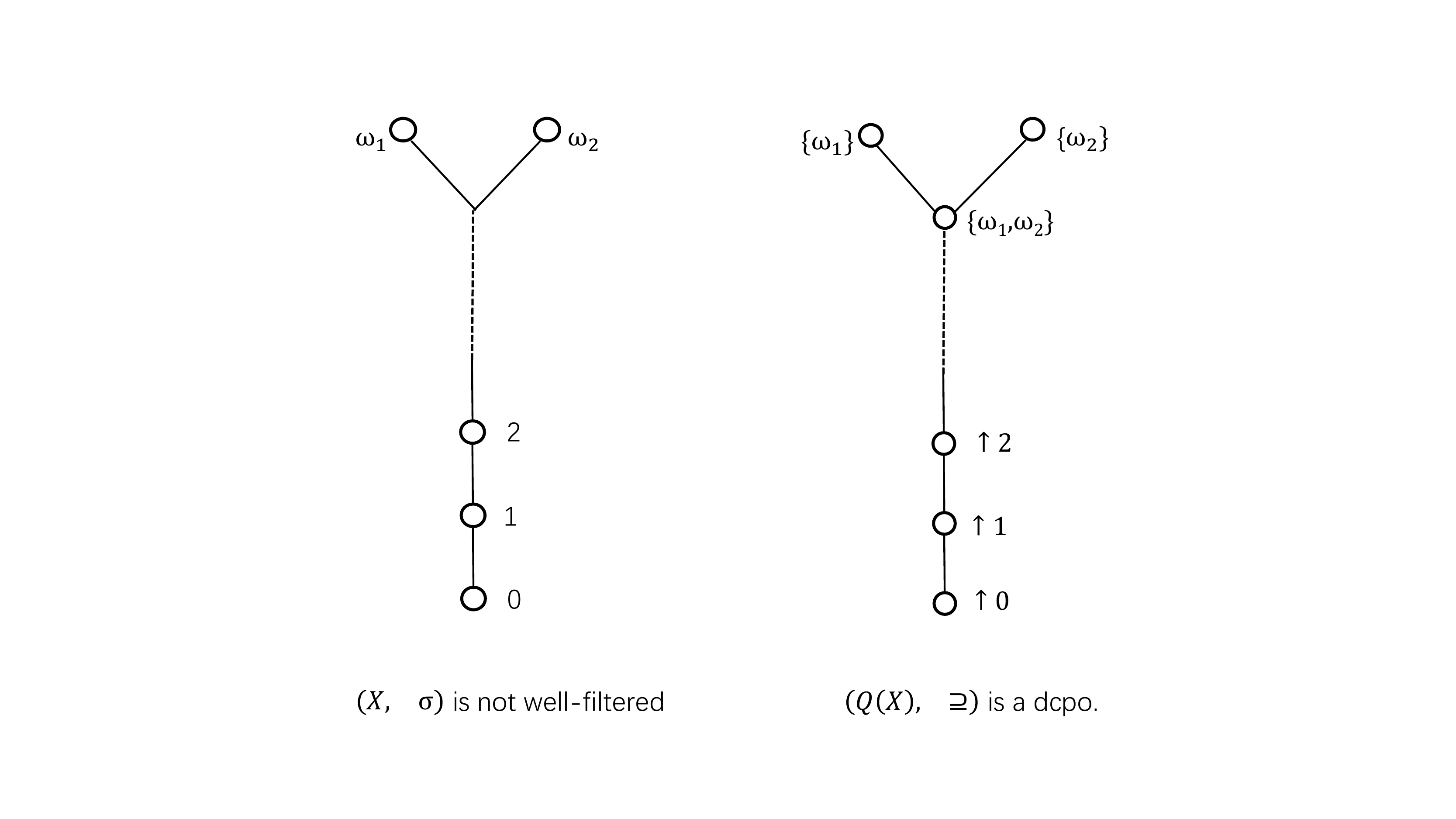}
\end{figure*}
\end{example}

\begin{theorem}\label{5}
 Let $X$ be a topological space. If $Q_{\upsilon}(X)$ is a $d^{*}$-space, then $X$ is a $d^{*}$-space.
 \end{theorem}
 \begin{proof}
 Suppose $D$ is a directed subset of $X$ with $\bigcap\limits_{d\in D}{\uparrow}d\cap{\uparrow}x\subseteq U$ for any $x\in X$ and any nonempty open set $U$ on $X$. Then $\mathcal{D}=\{{\uparrow}d\mid d\in D\}$ is directed in $Q(X)$. We claim that $\bigcap\limits_{d\in D} {\uparrow} _{Q_{\upsilon}(X)}({\uparrow}d )\cap{\uparrow}_{Q_{\upsilon}(X)}({\uparrow}x)\subseteq \Box U$. Indeed, suppose $K\in\bigcap\limits_{d\in D}{\uparrow} _{Q_{\upsilon}(X)}({\uparrow}d )\cap{\uparrow}_{Q_{\upsilon}(X)}({\uparrow}x)$. Then for all $d\in D$, ${\uparrow}d\leq K$ and ${\uparrow}x\leq K$. Thus for all  $ d\in D$,  $K\subseteq\uparrow d$ and $K\subseteq{\uparrow} x$. So $K\subseteq\bigcap\limits_{d\in D}{\uparrow} d\cap{\uparrow}x\subseteq U$, $K\in\Box U$. Since $U$ is nonempty, there exists $b\in U$ such that ${\uparrow}b\subseteq U$, which implies that $\Box U\neq\varnothing$. So there exists $d\in D$ such that ${\uparrow} _{Q_{\upsilon}(X)}({\uparrow} d )\cap{\uparrow}_{Q_{\upsilon}(X)}({\uparrow} x)\subseteq \Box U$. Assume $m\in{\uparrow}d\cap{\uparrow}x$,  then ${\uparrow}m\in{\uparrow} _{Q_{\upsilon}(X)}({\uparrow}d )\cap{\uparrow}_{Q_{\upsilon}(X)}({\uparrow}x)\subseteq \Box U$, it follows that ${\uparrow}m\subseteq U$. Hence, $m\in U$. We conclude that ${\uparrow}d\cap{\uparrow}x\subseteq U$.
 \end{proof}
 The following example demonstrates that the converse of  Theorem \ref{5} does not hold.
 \begin{example}
  Consider the set of natural numbers $N$ with the co-finite topology $\tau_{1}=\{U \mid N\backslash U$ is finite$\}\cup\{\varnothing\}$ and the single point set $\{a\}$ with the discrete topology $\tau_{2}=\{\varnothing,\{a\}\}$. The topology on $N\cup \{a\}$ is generated by the refinement of $\tau_{1}\vee\tau_{2}$. Obviously, it is a $T_{1}$ space because every singleton set is a closed set. Thus it is a $d^{*}$ space. Clearly, $Q(X)=2^{N\cup{\{a\}}}\backslash\{\varnothing\}$. Let $K_{0}=N$, $K_{1}=N\backslash\{0\}$, $K_{2}=N\backslash\{0,1\},\cdots, K_{n}=N\backslash\{0,1,2 \cdots n-1\}$,$\cdots$. Obviously, $\{K_{n}\}_{n\in N}$ is  directed in $Q(X)$. Let $K=K_{0}$. Then $\varnothing=\bigcap\limits_{n\in N}{\uparrow}_{Q_{\upsilon}(X)}K_{n}\cap{\uparrow}_{Q_{\upsilon}(X)}K=\bigcap\limits_{n\in N}{\uparrow}_{Q_{\upsilon}(X)}K_{n}=\varnothing\subseteq\Box\{a\}$. But for all $ n\in N$, $K_{n}\nsubseteq \{a\}$, that is, $K_{n}\notin\Box\{a\}$. Furthermore, for all $ n\in N$, we have ${\uparrow}_{Q_{\upsilon}(X)}K_{n}\cap{\uparrow}_{Q_{\upsilon}(X)}K={\uparrow}_{Q_{\upsilon}(X)}K_{n}\nsubseteq\Box\{a\}$. So $Q_{\upsilon}(X)$ is not a $d^{*}$-space.
  \end{example}
  The following example demonstrates that $Q_{\upsilon}(X)$ is a $d^{*}$-space but $X$ with the specialization order is  not a complete lattice.
 \begin{example}
 Consider the natural number $N$ with the co-countable topology $\tau_{c}=\{U\mid N\backslash U$ is countable $\}$. $(N,\tau_{c})$ is well filtered. We find that $Q_{\upsilon}(N)$ is a $d^{*}$-space but $N$ with the specialization order is  not a complete lattice. Indeed, suppose $\mathcal{K}$ is a directed subset of $Q(N)$ with $\bigcap_{K\in\mathcal{K}}{\uparrow_{Q_{\upsilon}(N)}K}\cap{\uparrow_{Q_{\upsilon}(N)}K^{'}}\subseteq\mathcal{U}$ for any $K^{'}$ in $Q(N)$ and nonempty open subset $\mathcal{U}$. Now let us do tow cases. If $K\cap K^{'}=\varnothing$ for some $K\in\mathcal{K}$, $\varnothing={\uparrow_{Q_{\upsilon}(K)}K}\cap{\uparrow_{Q_{\upsilon}(K)}K^{'}}\subseteq \mathcal{U}$. If $K\cap K^{'}\neq\varnothing$ for all $K\in\mathcal{K}$, $\bigcap_{K\in\mathcal{K}}{\uparrow_{Q_{\upsilon}(K)}K}\cap{\uparrow_{Q_{\upsilon}(K)}K^{'}}=\bigcap_{K\in\mathcal{K}}{\uparrow_{Q_{\upsilon}(K)}K\cap K^{'}}\subseteq\mathcal{U}$. Since $Q_{\upsilon}(N)$ is well filtered and $K\cap K^{'}$ is finite, there exists $K$ in $\mathcal{K}$ such that $ {\uparrow_{Q_{\upsilon}(K)}K\cap K^{'}}\subseteq\mathcal{U}$. So ${\uparrow_{Q_{\upsilon}(K)}K}\cap{\uparrow_{Q_{\upsilon}(K)}K^{'}}\subseteq \mathcal{U}$ for some $K$ in $\mathcal{K}$.
\end{example}

\bibliographystyle{./entics}

\end{document}